\documentclass[12pt]{article}

\usepackage{tikz}
\usepackage{listings}

\begin{document}
\title{Average angles of triangle in regular polygons}
\author{
        \textbf{Herman Muzychko\footnote{\textit{E-Mail: }herman.muzychko@tum.de}}
        \\\\
        Student at \\
        Physik-Department\\
        \textbf{Technische Universität München}\\
        James-Franck-Str. 1\\
        85748 Garching, Germany
}
\date{\today}

\maketitle

\begin{abstract}
In this article we will consider average angles of triangle, which share the same side with regular polygons. In particular we will count average angles in the triangle, which share the same bottom side with a square with length side $d=1$.
\end{abstract}

\section{Introduction}
Whenever we are talking about something average of continuity, we use integrals. 
\paragraph{Consider following problem:}
Let us count the average number of $\left\lbrace 1, 2, 3, 4, 5\right\rbrace$. We know, it will be $\frac{1+2+3+4+5}{5}=3$ — sum of numbers divided by their quantity. \\
But, we can solve this problem using integral. Let $x$ be variable, which represents each number and the limits of integral be from 1 to 5. Then, according to integrate rules we get:
$$\frac{1}{5-1}\int_{1}^{5} x dx=\frac{1}{4}\frac{x^{2}}{2}\Biggr|^{5}_{1}=\frac{25-1}{8}=3$$

\section{Problem}\label{previous work}
Consider a square in Cartesian coordinates and point $P(x, y)$. This point can freely move inside square. Position of $P(x, y)$ is connected to sides $a, b$ and angles $\alpha, \beta, \gamma$. Now using law of cosines and $c=1$ we get:

\tikzset{every picture/.style={line width=0.75pt}} %set default line width to 0.75pt        

\begin{tikzpicture}[x=0.75pt,y=0.75pt,yscale=-1,xscale=1]

%Shape: Square [id:dp9334721968730217] 
\draw   (41.5,28) -- (263,28) -- (263,249.5) -- (41.5,249.5) -- cycle ;
%Straight Lines [id:da02638533462944359] 
\draw    (96.5,95.5) -- (41.5,249.5) ;

%Straight Lines [id:da6628664967670366] 
\draw    (96.5,95.5) -- (263,249.5) ;

%Curve Lines [id:da30240872323683043] 
\draw    (88.5,118) .. controls (88.5,118) and (105.5,130) .. (116.5,113) ;

%Shape: Circle [id:dp053479021227356593] 
\draw  [fill={rgb, 255:red, 0; green, 0; blue, 0 }  ,fill opacity=1 ] (95.12,94.92) .. controls (95.44,94.16) and (96.32,93.8) .. (97.08,94.12) .. controls (97.84,94.44) and (98.2,95.32) .. (97.88,96.08) .. controls (97.56,96.84) and (96.68,97.2) .. (95.92,96.88) .. controls (95.16,96.56) and (94.8,95.68) .. (95.12,94.92) -- cycle ;
%Curve Lines [id:da4659228726429592] 
\draw    (51.5,225) .. controls (46.5,223) and (67.5,229) .. (69.5,249) ;

%Curve Lines [id:da6523134963232093] 
\draw    (226.5,249) .. controls (226.5,249) and (226.5,235) .. (242.5,230) ;

% Text Node
\draw (62,152) node   [align=left] {a};
% Text Node
\draw (175,148) node   [align=left] {b};
% Text Node
\draw (144,264) node   [align=left] {c};
% Text Node
\draw (39,263) node   [align=left] {(0,0)};
% Text Node
\draw (265,264) node   [align=left] {(1,0)};
% Text Node
\draw (39,15) node   [align=left] {(1,0)};
% Text Node
\draw (265,16) node   [align=left] {(1,1)};
% Text Node
\draw (105,130) node    {$\gamma $};
% Text Node
\draw (100,80) node    {$P( x,\ y)$};
% Text Node
\draw (400,81) node  [font=\large]  {$c^{
2} =a^{
2} +b^{
2} -2abcos( \gamma )$};
% Text Node
\draw (427,146) node  [font=\large]  {$ \begin{array}{l}
c\equiv 1\ \Rightarrow \alpha =\arccos\left[\frac{b^{2} +1-a^{2}}{2b}\right]\\
\end{array}$};
% Text Node
\draw (218,231) node    {$\alpha $};
% Text Node
\draw (75,228) node    {$\beta $};
% Text Node
\draw (400,51) node  [font=\large]  {$b^{
2} =a^{
2} +c^{
2} -2accos( \beta )$};
% Text Node
\draw (401,21) node  [font=\large]  {$a^{
2} =b^{
2} +c^{
2} -2bccos( \alpha )$};
% Text Node
\draw (469,188) node  [font=\large]  {$ \begin{array}{l}
\beta =\arccos\left[\frac{a^{2} -b^{2} +1}{2a}\right]\\
\end{array}$};
% Text Node
\draw (470,231) node  [font=\large]  {$ \begin{array}{l}
\gamma =\arccos\left[\frac{a^{2} +b^{2} -1}{2ab}\right]\\
\end{array}$};

\end{tikzpicture}

Using Cartesian coordinates, we define $a$ and $b$ sides by following:
\begin{center}
$a=\sqrt{x^{2}+y^{2}}$ and $b=\sqrt{(x-1)^{2}+y^{2}}$
\end{center}
At last inserting $a$ and $b$ in $\alpha, \beta, \gamma$ formuls we get:

\begin{equation}
\alpha=\arccos{\Biggl[\frac{x^{2}-x+y^{2}}{\sqrt{(x-1)^{2}+y^{2}}}\Biggl]}
\end{equation}
\begin{equation}
\beta=\arccos{\Biggl[\frac{x}{\sqrt{x^{2}+y^{2}}}\Biggl]}
\end{equation}
\begin{equation}
\gamma=\arccos{\Biggl[\frac{x^{2}-x+y^{2}}{\sqrt{x^{2}+y^{2}}\sqrt{(x-1)^{2}+y^{2}}}\Biggl]}
\end{equation}

Using formula (1) we can calculate exact angle $\alpha$ between $b$ and $c$ sides in every triangle with bottom side $c=1$. Similar comes for $\beta$ and $\gamma$. Now we are returning to our problem to calculate average $\alpha, \beta, \gamma$ angles inside $1\times1$ square.
\\
To find average angle we need to count sum of angles and divide by its count, on the other hand we can and will use another way -- integration, because it will bring us to exact answer. We use from introduction average integration method in regular polygon, hence we need double integral, in our case it is integrals along $x$ and $y$ sides of our square:
\begin{equation}
<\alpha>=I_{\alpha}=\int_{0}^{1}\int_{0}^{1} dy\,dx \arccos{\Biggl[\frac{x^{2}-x+y^{2}}{\sqrt{(x-1)^{2}+y^{2}}}\Biggl]}
\end{equation}
\begin{equation}
<\beta>=I_{\beta}=\int_{0}^{1}\int_{0}^{1} dy\,dx \arccos{\Biggl[\frac{x}{\sqrt{x^{2}+y^{2}}}\Biggl]}
\end{equation}
\begin{equation}
<\gamma>=I_{\gamma}=\int_{0}^{1}\int_{0}^{1} dy\,dx \arccos{\Biggl[\frac{x^{2}-x+y^{2}}{\sqrt{x^{2}+y^{2}}\sqrt{(x-1)^{2}+y^{2}}}\Biggl]}
\end{equation}
This integrals look very hard to calculate. We can simplify them by using polar coordinates $x=r\cos\varphi$, $y=r\cos\varphi$ and then expansion in 0 (Taylor series), but it is still a lot to do. Within this article we will not take this integral. Instead we will use computer algorithm.
\\
The easiest way to calculate numerically is to use computational machine. We will write a code on Python to deal with this problem.

\newpage
\section{Code (Python 3)}\label{code}
In this piece of code we are using the first method of counting average -- sum divided by quantity. More counts -- more accuracy. Thank to power of computers we can calculate thousands of decimal places.
\lstset{language=Python}
\begin{lstlisting}
from math import sqrt, acos, pi

# x from 0 to 1
# y from 0 to 1

sum_alpha = 0
sum_beta = 0
sum_gamma = 0
i = 0

for x in range(1, 1001):
    for y in range(1, 1001):
        a = sqrt((x/1000)**2 + (y/1000)**2)
        b = sqrt((x/1000 - 1)**2 + (y/1000)**2)
        
        alpha = acos((b**2 + 1 - a**2)/(2*b))*180/pi
        beta = acos((a**2 + 1 - b**2)/(2*a))*180/pi
        gamma = acos((a**2 + b**2 - 1)/(2*a*b))*180/pi
        
        sum_alpha += alpha
        sum_beta += beta
        sum_gamma += gamma
        
        i += 1
        
print("alpha =", sum_alpha/i)
print("beta =", sum_beta/i)
print("gamma =", sum_gamma/i)
\end{lstlisting}

\newpage
\section{Results}\label{results}
Using simple algorithm to count average angles, we have got these results:
\begin{center}
$<\alpha> = 45.064834706400624^{\circ}$\\
$<\beta> = 45.00000000000093^{\circ}$\\
$<\gamma> = 89.93516529359972^{\circ}$
\end{center}
It looks like, as if point $P(x, y)$ was in center of square and we measure angles of its triangle. 

\section{Conclusions}\label{conclusion}
The first thing comes to our minds is a principle of symmetry. If we consider geometric polygon, such as equilateral triangle, square, circle, etc., we get the same thing -- average angles will be only, as if $P(x, y)$ will be in geometric center of an object.
\\

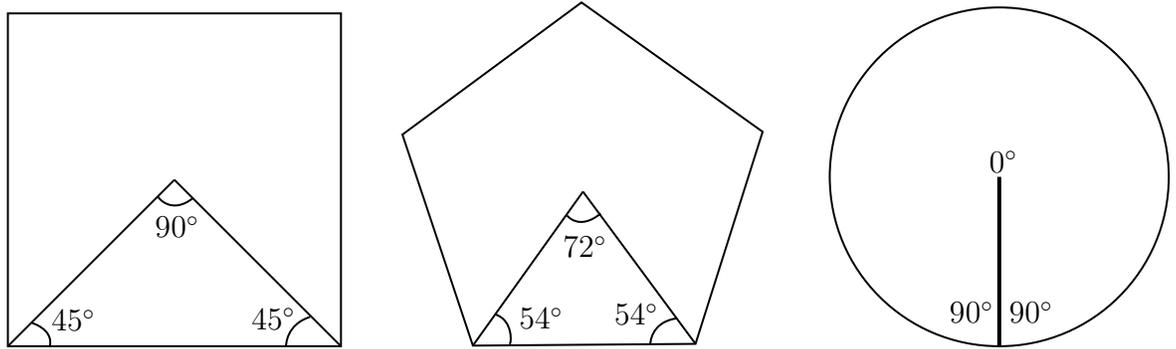
\tikzset{every picture/.style={line width=0.75pt}} %set default line width to 0.75pt        

\begin{tikzpicture}[x=0.75pt,y=0.75pt,yscale=-1,xscale=1]

%Shape: Square [id:dp9334721968730217] 
\draw   (11.5,23) -- (179.5,23) -- (179.5,191) -- (11.5,191) -- cycle ;
%Straight Lines [id:da02638533462944359] 
\draw    (95.5,107) -- (11.5,191) ;

%Straight Lines [id:da6628664967670366] 
\draw    (95.5,107) -- (179.5,191) ;

%Curve Lines [id:da30240872323683043] 
\draw    (87.35,115.53) .. controls (87.35,116.29) and (94.17,124.63) .. (104.79,116.29) ;

%Curve Lines [id:da4659228726429592] 
\draw    (23.64,179.24) .. controls (23.64,179.24) and (34.25,182.28) .. (32.74,190.62) ;

%Curve Lines [id:da6523134963232093] 
\draw    (151.82,190.62) .. controls (151.82,190.62) and (151.82,180) .. (163.95,176.21) ;

%Shape: Regular Polygon [id:dp7314276507271047] 
\draw   (358.31,189.81) -- (246.03,190.69) -- (210.49,84.18) -- (300.81,17.47) -- (392.16,82.75) -- cycle ;
%Straight Lines [id:da6511937353183006] 
\draw    (301.56,112.98) -- (246.03,190.69) ;

%Straight Lines [id:da6072143934600331] 
\draw    (301.56,112.98) -- (358.31,189.81) ;

%Curve Lines [id:da6188948899531554] 
\draw    (257.5,175) .. controls (266.5,178) and (265.5,190) .. (264.5,190) ;

%Curve Lines [id:da8655069263177193] 
\draw    (293.5,125) .. controls (293.5,125) and (299.5,133) .. (310.5,124) ;

%Curve Lines [id:da03756264287596012] 
\draw    (335.5,190) .. controls (335.5,190) and (335.5,179) .. (348.5,177) ;

%Shape: Circle [id:dp42841542746400973] 
\draw   (426,105.5) .. controls (426,58.28) and (464.28,20) .. (511.5,20) .. controls (558.72,20) and (597,58.28) .. (597,105.5) .. controls (597,152.72) and (558.72,191) .. (511.5,191) .. controls (464.28,191) and (426,152.72) .. (426,105.5) -- cycle ;
%Straight Lines [id:da13159240816346984] 
\draw [line width=1.5]    (511.5,105.5) -- (511.5,191) ;

% Text Node
\draw (96.83,130.7) node    {$90^{\circ} $};
% Text Node
\draw (145.37,176.97) node    {$45^{\circ} $};
% Text Node
\draw (44.49,177.73) node    {$45^{\circ} $};
% Text Node
\draw (280.49,176.73) node    {$54^{\circ} $};
% Text Node
\draw (328.49,174.73) node    {$54^{\circ} $};
% Text Node
\draw (302.49,140.73) node    {$72^{\circ} $};
% Text Node
\draw (497.49,173.73) node    {$90^{\circ} $};
% Text Node
\draw (527.49,173.73) node    {$90^{\circ} $};
% Text Node
\draw (513.49,97.73) node    {$0^{\circ} $};

\end{tikzpicture}

Therefore we can count every average angle in triangle in every regular polygon using empirical formula in degrees:
\begin{equation}
<\alpha> = \frac{(n-2)}{2n}\cdot180^{\circ}
\end{equation}
\begin{equation}
<\beta> = \frac{(n-2)}{2n}\cdot180^{\circ}
\end{equation}
\begin{equation}
<\gamma> = 180^{\circ} - <\alpha> - <\beta>
\end{equation}
\textbf{Note:}
$n$ -- number of sides. For circle $n\rightarrow\infty$
\\
Lastly we can use (7), (8) and (9) for vector problems in mechanics, for instance for momentum of particles after collisions, etc. 

\bibliographystyle{abbrv}
\bibliography{main}

\end{document}